\documentclass[12 pt, reqno]{amsart}

\usepackage{amsthm,amssymb,amstext,amscd,amsfonts,amsbsy,amsxtra,latexsym,cite}
\usepackage{mathrsfs}
\usepackage{verbatim}
\usepackage[latin1]{inputenc}
\usepackage{graphicx, graphics, wrapfig}

\newcommand{\field}[1]{\mathbb{#1}}

\newcommand{\N}{\field{N}}
\newcommand{\Z}{\field{Z}}

\newcommand{\C}{\field{C}}
\def\({\left(}
\def\){\right)}

\newcommand{\sC}{\mathscr{C}}

\newcommand{\qbinom}[2]{\genfrac{[}{]}{0pt}{}{#1}{#2}}

\theoremstyle{plain}
\newtheorem{theorem}{Theorem}
\newtheorem*{theorem*}{Theorem}
\newtheorem{lemma}[theorem]{Lemma}

\newtheorem*{conjecture*}{Conjecture}

\theoremstyle{definition}

\theoremstyle{remark}

\newtheorem*{remark*}{Remark}
\newtheorem*{remarks}{Remarks}

\numberwithin{theorem}{section} \numberwithin{equation}{section}

\begin{document}

\title{False theta functions and companions to Capparelli's identities}
\author{Kathrin Bringmann}
\address{Mathematical Institute\\University of
Cologne\\ Weyertal 86-90 \\ 50931 Cologne \\Germany}
\email{kbringma@math.uni-koeln.de}
\author{Karl Mahlburg}
\address{Department of Mathematics \\
Louisiana State University \\
Baton Rouge, LA 70802\\ U.S.A.}
\email{mahlburg@math.lsu.edu}


\date{\today}


\thanks{The research of the first author was supported by the Alfried Krupp Prize for Young University Teachers of the Krupp Foundation and the research leading to these results has received funding from the European Research Council under the European Union's Seventh Framework Programme (FP/2007-2013) / ERC Grant agreement n. 335220 - AQSER.  The second author was supported by NSF Grant DMS-1201435.}

\begin{abstract}
Capparelli conjectured two modular identities for partitions whose parts satisfy certain gap conditions, which were motivated by the calculation of characters for the standard modules of certain affine Lie algebras and by vertex operator theory.
These identities were subsequently proved and refined by Andrews,
who related them to Jacobi theta functions, and also by Alladi-Andrews-Gordon, Capparelli and Tamba-Xie.
In this paper we prove two new companions to Capparelli's identities, where the evaluations are expressed in terms of Jacobi theta functions and false theta functions.
\end{abstract}

\maketitle

\section{Introduction and statement of results}

The study of hypergeometric $q$-series identities for partitions whose adjacent parts satisfy minimum gap conditions has a long and rich history, including Euler's Theorem for partitions into distinct parts (see e.g. Corollary 1.2 in \cite{And98}), and the famous Rogers-Ramanujan identities \cite{RR19}. The basic combinatorial problem of evaluating the resulting generating functions has deep ramifications in the theory of hypergeometric $q$-series and modular forms (cf. \cite{And74,Gor61}),
 as the resulting identities provide infinite products that are essentially (up to rational $q$-powers) meromorphic quotients of elliptic theta functions.

The partition identities discussed in this paper are inspired by Lepowsky and Wilson's seminal work in \cite{LW78}, which introduced vertex operators as a method for explicitly constructing affine Lie algebras, and \cite{LW81}, which extended this construction to $Z$-algebras, which are certain
generalized vertex-algebraic structures; see \cite{FLM88} for further developments.  From this perspective the construction of the corresponding standard modules is naturally linked to combinatorial partition identities, such as the new interpretation of the Rogers-Ramanujan identities proven in \cite{LW81}. Indeed, Lepowsky and Wilson's construction of the standard modules of $A^{(1)}_1$ results in formulas that coincide with the generalized Rogers-Ramanujan identities due to Andrews \cite{And74}, Bressoud \cite{Bre80}, and Gordon \cite{Gor61}.

Lepowsky and Milne \cite{LM78} also showed that the Rogers-Ramanujan identities arise in character formulas for the level 2 standard modules for $A^{(2)}_2$, which was later proven using $Z$-algebras by Capparelli \cite{Cap93}.
Capparelli additionally used $Z$-algebras to conjecturally construct the level 3 standard modules (also see \cite{Cap88}), and found two striking formulas for the generating functions of partitions satisfying certain gap conditions and smallest part restrictions.
Capparelli's identities were a significant development in the theory of vertex operator algebras, as they were the first notable examples that had not been previously appeared in the literature of partition identities, but were instead discovered using vertex-operator-theoretic techniques. Independent proofs of Capparelli's conjectures were subsequently given by Andrews, Alladi-Andrews-Gordon, Capparelli and Tamba-Xie, as we discuss below; the fact that the identities intersect with a wide variety of fields further indicates their significance. Many other such identities have since also been found; for example, see [18, 19, 23] for further discussion of the role of affine Lie algebras, vertex operator methods, and statistical mechanics.

Following the  then-conjectural statement of Theorem 21 of \cite{Cap93}, we say that a partition satisfies Capparelli's {\it level $3$ gap condition} if successive parts differ by at least $2$, and two parts differ by $2$ or $3$ only if their sum is a multiple of $3$.
The level $3$ gap condition can also be equivalently written in terms of part multiplicities, which is one form of the combinatorial identities that naturally occur in the study of character formulas. For example, the following formulation is the special case $k=1$ of (11.2.8) in \cite{MP99}, which describes the partition ideals that arise from root lattices. Define indicator functions such that $\psi_j(\lambda) = 1$ if $j$ is a part of $\lambda$, and $\psi_j(\lambda) = 0$ otherwise. The level $3$ gap condition is satisfied if and only if for all $j \geq 1$,
\begin{align*}
\psi_{3j+2} + \psi_{3j} + \psi_{3j-1} & \leq 1, \\
\psi_{3j+1} + \psi_{3j} + \psi_{3j-2} & \leq 1, \\
\psi_{3j-1} + \psi_{3j-2} & \leq 1.
\end{align*}
Note that (11.2.8) of \cite{MP99} is actually a system of four inequalities, but in the special case $k=1$ it is overdetermined and reduces to the above.

In order to state Capparelli's conjectures, now identities, we also require enumeration functions for the partitions described above. Let $c_m(n)$ denote the number of partitions of $n$ that satisfy the level $3$ gap condition and whose parts are all at least $m$.
Capparelli also considered the closely related function
$c^\ast_2(n)$, which denotes the number of partitions that satisfy the level $3$ gap conditions and additionally do not contain $2$ as a part.
Note that this can be expressed in terms of the $c_m(n)$ through the simple combinatorial relation
\begin{equation}
\label{E:c2*}
c^\ast_2(n) = c_1(n) - c_2(n) + c_3(n).
\end{equation}

Capparelli's identities are stated in the following theorem, which was proven using techniques from $q$-series by Andrews \cite{And94} (this paper contains only the first of the two identities), and Andrews, Alladi, and Gordon \cite{AAG95}. The identities were also proven using the $Z$-algebra program of \cite{LW81} by both Tamba and Xie \cite{TX95} and Capparelli himself \cite{Cap96}.
For $j = 1,2$, let $d_j(n)$, denote the number of partitions of $n$ into distinct parts that are not $\pm j \pmod{6}$.
\begin{theorem*}[\cite{AAG95, And94,Cap96,TX95}]
For $n \geq 1$,
\begin{align}
c_2(n) & = d_1(n), \label{E:c2=d1} \\
c^\ast_2(n) & = d_2(n). \label{E:c2*=d2}
\end{align}
\end{theorem*}
\begin{remark*}
The above theorem is not the original combinatorial formulation of Capparelli's identities, but it is an elementary exercise in infinite product generating functions to verify that, for example, $d_1(n)$ also enumerates the number of partitions of $n$ into parts congruent to $\pm 2, \pm 3 \pmod{12},$ as in Theorem 21 A of \cite{Cap93}.
\end{remark*}

However, despite their natural connections to Lie theory, the study of Capparelli's identities using hypergeometric $q$-series gives additional information that is of particular number-theoretic interest.
The preceding combinatorial results can be equivalently stated as generating function identities, and for $m \geq 1$ we write
\begin{equation*}
\sC_m(q) := \sum_{n \geq 0} c_m(n) q^n = \sum_{\substack{\lambda \text{ level } 3 \text{ gaps} \\ \lambda_j \geq m}} q^{|\lambda|}.
\end{equation*}
Here $|\lambda|$ denotes the {\it size} of a partition $\lambda$, and the summation subscript in the final expression is an abbreviation for ``$\lambda$ satisfying the level $3$ gap condition, with all parts of $\lambda$ of at least size $m$.''
Similarly, we set
\begin{equation*}
\sC^\ast_2(q) := \sum_{n \geq 0} c^\ast_2(n)q^n.
\end{equation*}
Note that \eqref{E:c2*} is equivalently expressed as
\begin{equation}
\label{E:C2*q}
\sC^\ast_2(q) = \sC_1(q) - \sC_2(q) + \sC_3(q).
\end{equation}

In this paper, we adopt the standard $q$-factorial notation for $a \in \C$ and $n\in\N_0 \cup\{\infty\}$, namely  $(a)_n = (a;q)_n := \prod_{j = 0}^{n-1} (1-aq^j).$  We also use the additional shorthand $(a_1, \dots , a_r)_n := (a_1)_n \cdots (a_r)_n$. We can now rewrite \eqref{E:c2=d1} and \eqref{E:c2*=d2} as
\begin{align}
\sC_2(q) & = \left(-q^2, -q^3, -q^4, -q^6; q^6\right)_\infty, \label{E:C2=prod} \\
\sC^\ast_2(q) & = \left(-q, -q^3, -q^5, -q^6; q^6\right)_\infty. \label{E:C2*=prod}
\end{align}
These are {\it modular} identities in the sense that the right-hand sides are essentially weakly holomorphic modular forms.

In fact, Andrews' results in \cite{And94} also include a refinement of \eqref{E:C2=prod} that is of additional number-theoretic interest due to the presence of an additional parameter.  Moreover, in \cite{AAG95} Alladi, Andrews, and Gordon proved a similar result for \eqref{E:C2*=prod}, as well as further refinements of both identities that provide additional combinatorial information. The general results in \cite{AAG95} are best stated as identities for {\it three-colored} partitions, in which each part may be labeled with one of three distinct colors; identities for three-colored partitions also arise in \cite{MP01}, where the basic $A_2^{(1)}$-module is constructed using vertex operator methods
We further note that in \cite{Sil04} Sills proved a one-parameter generalization of an ``analytic counterpart'' to Capparelli's identities, using Bailey chains to obtain interesting hypergeometric $q$-series representations for infinite products related to \eqref{E:C2=prod} -- \eqref{E:C2*=prod}.

However, in our present study we focus on the one-parameter refinements of Capparelli's identities, as we are particularly interested in the automorphic properties of $q$-series.
In order to describe the refined identities, for $j \in \{1, 2\}$, we let $\nu_j(\lambda)$ be the number of parts of $\lambda$ that are congruent to $j$ modulo $3$. The refined Capparelli generating functions are then defined as
\begin{equation*}
\sC_m(t;q) := \sum_{\substack{\lambda \text{ level } 3 \text{ gaps} \\ \lambda_j \geq m}} t^{\nu_1(\lambda) - \nu_2(\lambda)} q^{|\lambda|},
\end{equation*}
and similarly the refined generating function for partitions without $2$ as a part is
\begin{equation*}
\sC^\ast_2(t;q) := \sC_1(t;q) - \sC_2(t;q) + \sC_3(t;q).
\end{equation*}

The refinements of Capparelli's identities provide product identities for these generating functions. We write the results in terms of the {\it Jacobi theta function}
\begin{equation}
\label{E:theta}
\theta(z;q) := \left(-z, -z^{-1}q, q; q\right)_\infty =
\sum_{k \in \Z} z^k q^{\frac{k(k-1)}{2}},
\end{equation}
where the final equality follows from Jacobi's Triple Product identity ((2.2.10) in \cite{And98}).
\begin{theorem*}[\cite{And94,AAG95}]
The following identities hold:
\begin{align}
\sC_2(t;q) &=  \frac{\theta\left(tq^4; q^6\right)}{\left(q^3; q^3\right)_\infty}
\label{E:C2tq=prod}, \\
\sC^\ast_2(t;q) & = \frac{\theta\left(tq; q^6\right)} {\left(q^3; q^3\right)_\infty}
\label{E:C2*tq=prod}.
\end{align}
\end{theorem*}
\begin{remark*}
Although the first of these identities is not stated in \cite{And94}, it follows implicitly as a limiting case of Theorem 2 from that paper. It is also a special case of Theorem 1 from \cite{AAG95} (specifically, see Remark 1 on page 646 of \cite{AAG95} and set $a = t, b = t^{-1}$).  The second identity is not stated in \cite{AAG95}, but is implied by Theorem 3 and the discussion in the Acknowledgments Section on page 658.
\end{remark*}
\begin{remark*}
Building on the previous discussion of modularity, we note that $\theta(z;q)$ is essentially a holomorphic Jacobi form as introduced by Eichler and Zagier \cite{EZ85}.
\end{remark*}

Our main result provides two new evaluations for the generating functions of partitions satisfying the level $3$ gap conditions as introduced by Capparelli.

\begin{theorem}
\label{T:Main}
Let $\chi_3$ denote the shifted Dirichlet character defined by $\chi_3(m) := \left(\frac{m+1}{3}\right).$ Then following identities hold:
\begin{align}
\sC_1(t;q) & = \left(-q^3; q^3\right)_\infty \theta\left(-t^2 q^2; q^6\right)
+ \sC_2(t;q) \big(1 - \Theta_1(t;q)\big) + \sC^\ast_2(t;q) \big(1 - \Theta_2(t;q)\big),
\nonumber
\\
\sC_3(t;q)  & =
-\left(-q^3; q^3\right)_\infty \theta\left(-t^2q^2; q^6\right)
+ \sC_2(t;q) \Theta_1(t;q) + \sC^\ast_2(t;q) \Theta_2(t;q), \label{E:C3tq=prod}
\end{align}
where
\begin{align*}
\Theta_1(t;q) & := \sum_{k \geq 0} \chi_3(k) t^{-k} q^{k(k+1)}, \\
\Theta_2(t;q) & := \sum_{k \geq 0} \chi_3(k) t^k q^{k^2}.
\end{align*}
\end{theorem}
\begin{remarks}
{\it 1.}  As described in an earlier remark, the factor $\theta\left(-t^2q^2; q^6\right)$ is essentially a Jacobi form. The $\Theta_j$ should be thought of as ``false'' analogues of Jacobi's theta function, in the sense that they are given as half-lattice sums that are twisted by shifted characters. Note that the first series, for example, is equivalently written as
\begin{equation*}
\Theta_1(t;q) := \sum_{k \geq 0}\Big(t^{-3k} q^{3k(3k+1)} - t^{-(3k+1)} q^{(3k+1)(3k+2)}\Big).
\end{equation*}
Following Zagier's introduction of {\it quantum modular forms} in the seminal paper \cite{Zag10}, there has been a great deal of recent work illuminating the connections between false theta functions and classical automorphic forms; for example, see \cite{FOR13}. \\

\noindent {\it 2.} When $t = 1$, these formulas are reminiscent of the ``bosonic'' evaluations that are frequently found in the study of solvable lattice models and/or characters for affine Lie algebras; for example, see Section 2 of \cite{War96}. Furthermore, the appearance of false theta functions is not entirely unexpected; see \cite{BMO96} for other examples arising from mathematical physics, as well as \cite{MS08} for examples in $q$-series derived from Bailey pairs. However, the ``Jacobi-like'' parameter $t$ is often absent in identities that are derived from Lie-theoretic or physical settings, and our formulas therefore provide additional information.
 \\

\noindent {\it 3.} In light of the identities \eqref{E:C2*q}, \eqref{E:C2tq=prod}, and \eqref{E:C2*tq=prod}, either one of our new formulas implies the other. In fact, it turns out to be convenient to prove the four formulas \eqref{E:C2tq=prod} -- \eqref{E:C3tq=prod} simultaneously.
\end{remarks}

The remainder of the paper is structured as follows. In the next section we recall several standard results from the theory of hypergeometric $q$-series. We then prove Theorem \ref{T:Main} in Section \ref{S:proof} by applying techniques from the theory of $q$-difference equations \cite{And75, GR90}.

\section*{Acknowledgments}

The authors thank Jim Lepowsky and Andrew Sills for providing several additional references and historical context. The authors also thank the two anonymous referees for their detailed remarks that significantly improved the exposition.

\section{Hypergeometric $q$-series identities}
\label{S:qseries}
In this section we record a number of identities that are useful in the evaluation of the generating functions that are the main topic of the paper. If $0 \leq m \leq n$, the {\it $q$-binomial coefficient} is denoted by
\begin{equation*}
\qbinom{n}{m}_q := \frac{(q;q)_n}{(q;q)_m (q;q)_{n-m}}.
\end{equation*}
We also need the limiting case
\begin{equation}
\label{E:qbinomLim}
\lim_{n \to \infty} \qbinom{n}{m}_q = \frac{1}{(q;q)_m}.
\end{equation}

Next, we recall two identities due to Euler, which state (see (2.2.5) and (2.2.6) in \cite{And98})
\begin{align}
\label{E:Euler1}
\frac{1}{(x;q)_\infty} & = \sum_{n \geq 0} \frac{x^n}{(q;q)_n}, \\
\label{E:Euler2}
(x;q)_\infty &= \sum_{n \geq 0} \frac{(-1)^n x^n q^{\frac{n(n-1)}{2}}}{(q;q)_n}.
\end{align}
A related identity is
\begin{equation}
\label{E:CauchyEuler1}
\sum_{\substack{n \geq 0 \\ n \text{ even}}} \frac{q^{\frac{n(n-1)}{2}}}{\left(q; q\right)_n}
= \frac{1}{\left(q; q^2\right)_\infty} = \left(-q; q\right)_\infty;
\end{equation}
the first equality follows from Cauchy's identity, which is (2.2.8) in \cite{And98}.

We also need the following identity from Ramanujan's famous ``Lost Notebook'', which appears as (4.1) in \cite{And79}:
\begin{equation}
\label{E:Ramid}
\sum_{n \geq 0} \frac{q^n}{(-aq;q)_n (-bq;q)_n}
= \left(1 + a^{-1}\right) \sum_{n \geq 0} \frac{(-1)^n q^{\frac{n(n+1)}{2}}\left(\frac{b}{a}\right)^n}{(-bq;q)_n}
- \frac{a^{-1} \sum_{n \geq 0} (-1)^n q^{\frac{n(n+1)}{2}}\left(\frac{b}{a}\right)^n}{(-aq, -bq; q)_\infty}.
\end{equation}
Finally, in order to derive expressions involving false theta functions, we use a related identity of Rogers \cite{Rog17}, which states that
\begin{equation}
\label{E:Rogid}
\sum_{n \geq 0}\frac{(-1)^n y^{2n} q^{\frac{n(n+1)}{2}}}{(yq;q)_n}
= \sum_{n \geq 0} (-1)^n y^{3n} q^{\frac{n(3n+1)}{2}} \left(1 - y^2 q^{2n+1}\right).
\end{equation}

\section{Finite evaluations and proof of Theorem \ref{T:Main}}
\label{S:proof}

In order to simultaneously work with $\sC_1, \sC_2, \sC_3$, and $\sC^\ast_2$, we incorporate the part indicator functions into the generating series.  We control the smallest part conditions through two parameters $\alpha, \beta \in \{0, 1\}$, which determine whether or not the indicators are ``active''. Specifically, if $\alpha = 0$, then no parts of size $1$ are allowed, but when $\alpha = 1$ such a part is permitted, and $\beta$ has a similar affect on parts of size $2$. The generating functions with these indicators taken into account are then defined by
\begin{equation*}
\sC^{\alpha, \beta}(t;q) := \sum_{\lambda \text{ level } 3 \text{ gaps}}
\big(1 - (1-\alpha)\psi_1(\lambda)\big) \big(1 - (1-\beta)\psi_2(\lambda)\big)
t^{\nu_1(\lambda) - \nu_2(\lambda)} q^{|\lambda|}.
\end{equation*}
This notation corresponds to the four series of interest as follows:
\begin{equation*}
\sC_1 = \sC^{1,1}, \quad\quad \sC_2 = \sC^{0,1}, \quad\quad \sC^\ast_2 = \sC^{1,0}, \quad\quad \sC_3 = \sC^{0,0}.
\end{equation*}

The results in Theorem \ref{T:Main} are equivalent to the following generating function evaluation.
\begin{theorem}
If $\alpha, \beta \in \{0, 1\}$, then
\label{T:Mainab}
\begin{align*}
\sC^{\alpha, \beta}(t;q)  =&
(\alpha + \beta - 1) \left(-q^3; q^3\right)_\infty \theta\left(-t^2 q^2; q^6\right)
+ \frac{\theta\left(tq^4; q^6\right)}{\left(q^3; q^3\right)_\infty} \Big(\beta + (1 - \alpha - \beta)
\Theta_1(t;q)\Big)
\notag \\
  &+ \frac{\theta\left(tq; q^6\right)}{\left(q^3; q^3\right)_\infty}\Big(\alpha + (1 - \alpha - \beta)\Theta_2(t;q)\Big).
\end{align*}
\end{theorem}
\begin{remark*}
In the two cases considered by Capparelli the above expression simplifies drastically, as $1 - \alpha - \beta = 0$.
\end{remark*}

We prove our generating function evaluations by following Andrews' approach in \cite{And94}, as well as arguments from \cite{AAG95}. In particular, we consider recurrences for finite truncations of the generating functions, defining for $M\in\N$
\begin{equation*}
C^{\alpha,\beta}_M(t;q) := \sum_{\substack{\lambda \text{ level } 3 \text{ gaps} \\ \lambda_j \leq M}}
\big(1 - (1-\alpha)\psi_1(\lambda)\big) \big(1 - (1-\beta)\psi_2(\lambda)\big)
t^{\nu_1(\lambda) - \nu_2(\lambda)} q^{|\lambda|}.
\end{equation*}
For convenience of notation, we regularly suppress the superscript, and sometimes the arguments, writing only $C_M(t;q)$ or $C_M$, with $\alpha$ and $\beta$ unspecified.

Theorem \ref{T:Mainab} arises as the limiting case $n \to \infty$ of a finite evaluation; we will calculate this limit in Section \ref{S:proof:infinite}.
\begin{lemma}
\label{L:Mainabfinite}
For $n \geq 0$ we have
\begin{align*}
& C^{\alpha,\beta}_{3n-2}(t;q) \\
& =\sum_{\substack{0 \leq j \leq n \\ j \equiv n \pmod{2}}}
q^{\frac{3j(j-1)}{2}}\qbinom{n}{j}_{q^3} \left(-t^{-1}q^2, -tq^4; q^6\right)_{\frac{n-j}{2}}
\left(\beta + t(1 - \alpha - \beta)\sum_{\ell=1}^{\frac{n-j}{2}} \frac{q^{6\ell-2}}{\left(-t^{-1}q^2, -tq^4; q^6\right)_\ell}\right) \notag \\
& \qquad + \sum_{\substack{0 \leq j \leq n - 1 \\ j \equiv n-1 \pmod{2}}}
q^{\frac{3j(j-1)}{2}}\qbinom{n}{j}_{q^3} \left(-t^{-1}q^5, -tq^7; q^6\right)_{\frac{n-1-j}{2}} \\
& \qquad \qquad \times \left(1 - \beta + \alpha tq - t(1 - \alpha - \beta)\sum_{\ell=1}^{\frac{n-1-j}{2}} \frac{q^{6\ell+1}}{\left(-t^{-1}q^5, -tq^7; q^6\right)_\ell}\right).
\end{align*}
\end{lemma}

\subsection{Proof of Lemma \ref{L:Mainabfinite}}
\label{S:proof:finite}

We prove Lemma \ref{L:Mainabfinite} by following the general framework of \cite{AAG95} and \cite{And94}, although the presence of the indicators $(\alpha, \beta)$ introduces intricate boundary effects throughout the calculations.
As described in equations (4.2) -- (4.4) of \cite{And94}, the finite generating functions satisfy the recurrences (for $n \geq 2$)
\begin{align}
\label{E:Cnrec}
C_{3n-1}(t;q) &= C_{3n-2}(t;q) + t^{-1} q^{3n-1} C_{3n-5}(t;q), \\
C_{3n}(t;q) &= C_{3n-1}(t;q) + q^{3n}C_{3n-3}(t;q), \notag \\
C_{3n+1}(t;q) &= C_{3n}(t;q) + tq^{3n+1}C_{3n-3}(t;q) + q^{6n}C_{3n-5}(t;q). \notag
\end{align}
Note that these recurrences arise from conditioning on the largest parts in a partition and therefore depend only on parts of size $5$ or larger. As the indicators $\alpha$ and $\beta$ only affect parts of size $1$ and $2$, Andrews recurrences' are thus unchanged. However, the indicators do affect the initial values, which are obtained by explicitly listing the partitions with largest part at most $4$ that satisfy the level $3$ gap condition:
\begin{align}
\label{E:Cninit}
C_1(t;q) & = 1 + \alpha t q, \\
C_2(t;q) & = 1 + \alpha t q + \beta t^{-1}q^2, \notag \\
C_3(t;q) & = 1 + \alpha t q + \beta t^{-1} q^2 + q^3, \notag \\
C_4(t;q) & = 1 + \alpha t q + \beta t^{-1} q^2 + q^3 + tq^4 + \beta q^6. \notag
\end{align}

Andrews also showed that when the recurrences \eqref{E:Cnrec} are iterated, they combine to imply the single recurrence
\begin{align}
\label{E:C3n+1rec}
C_{3n+1} = & \left(1 + q^{3n}\right) C_{3n-2} + \left(t^{-1} q^{3n-1} + tq^{3n+1} + q^{6n}\right)C_{3n-5}  + q^{6n-3}\left(1 - q^{3n-3}\right)C_{3n-8},
\end{align}
which holds for $n \geq 3.$ In fact, it is convenient to assume that this recurrence holds for $n \geq 2$, which is achieved by setting $n=1$ in the first line of \eqref{E:Cnrec} to obtain the value $C_{-2}(t;q) = \beta$. Combined with the known expressions for $C_1$ and $C_4$ in \eqref{E:Cninit}, these are the necessary initial values.

Our goal is to gather the $C_M$ in a series and thereby derive a solvable $q$-difference equation, but this has a simpler shape if we first renormalize by defining
\begin{equation}
\label{E:gammadef}
\gamma_n=\gamma_{n}(t;q) := \frac{C_{3n-2}(t;q)}{(q^3; q^3)_n}.
\end{equation}
After shifting $n \mapsto n-1$, the recurrence \eqref{E:C3n+1rec} then becomes
\begin{equation}
\label{E:gammarec}
\begin{split}
&\left(1-q^{3n}\right) \left(1 - q^{3n-3}\right) \gamma_n\\
&= \left(1 - q^{6n-6}\right) \gamma_{n-1}
+ \left(t^{-1} q^{3n-4} + t q^{3n-2} + q^{6n-6}\right) \gamma_{n-2}
+ q^{6n-9} \gamma_{n-3},
\end{split}\end{equation}
which now holds for $n \geq 3$. From \eqref{E:Cninit} we see that the initial values are
\begin{equation}
\label{E:gammainit}
\gamma_0 = \beta, \qquad \gamma_1 = \frac{1 + \alpha t q}{1 - q^3}, \qquad
\gamma_2 = \frac{1 + \alpha t q + \beta t^{-1} q^2 + q^3 + tq^4 + \beta q^6}{\left(1 -q^3\right) \left(1 - q^6\right)}.
\end{equation}

Now we set
\begin{equation}
\label{E:Fdef}
F(z) = F(z,t;q) := \sum_{n \geq 0} \gamma_n(t;q) z^n,
\end{equation}
and we also briefly use the notation
\begin{equation*}
F^{(m)}(z) := \sum_{n \geq m} \gamma_n(t;q) z^n
\end{equation*}
for the series with truncated initial terms.
Multiplying \eqref{E:gammarec} by $z^n$ and summing over $n \geq 3$, we obtain the series identity
\begin{align*}
F^{(3)}(z) - \left(1 + q^{-3}\right)& F^{(3)}\left(zq^3\right) + q^{-3}F^{(3)}\left(zq^6\right) =
zF^{(2)}(z) - zF^{(2)}\left(zq^6\right) \\
& + z^2\left(t^{-1}q^2 + tq^4\right) F^{(1)}\left(zq^3\right) + z^2 q^6 F^{(1)}\left(zq^6\right) + z^3 q^9 F\left(zq^6\right).
\end{align*}
Using the initial values in \eqref{E:gammainit} to add back in the missing coefficients of $z^0, z^1,$ and $z^2$, we obtain the $q$-difference equation
\begin{align}
(1-z)F(z) = & \left(1 + q^{-3} + z^2 t^{-1} q^2 + z^2 t q^4\right) F\left(zq^3\right)
+ \left(1 + zq^3\right)\left(-q^{-3} + z^2 q^6\right) F\left(zq^6\right) \notag \\
& \qquad \qquad \qquad + z^2 (1 - \alpha - \beta) t q^4.
\label{E:Fqdiff}
\end{align}
Note that while the $z^0$ and $z^1$ terms vanish in all cases, the $z^2$ term only vanishes in the two cases considered by Capparelli (recall \eqref{E:C2tq=prod} and \eqref{E:C2*tq=prod})).

We renormalize once more by setting
\begin{equation}
\label{E:Hdef}
H(z) := \frac{F(z)}{(-z;q^3)_\infty},
\end{equation}
and obtain
\begin{align}
\label{E:Hqdiff}
\left(1 - z^2\right)H(z) = & \left(1 + q^{-3} + z^2 t^{-1} q^2 + z^2 tq^4\right) H\left(zq^3\right)
+ \left(-q^{-3} + z^2 q^6\right) H\left(zq^6\right) \notag \\
& \qquad \qquad + (1 - \alpha - \beta) t \sum_{n \geq 0} \frac{(-1)^n z^{n+2} q^{3n+4}}{\left(q^3; q^3\right)_n}.
\end{align}
For the last term we expanded $(-zq^3; q^3)^{-1}_\infty$ as a series using \eqref{E:Euler1}.

At this point our approach varies from \cite{AAG95}, where a solution to \eqref{E:Fqdiff} was found by using the general theory of second order $q$-difference equations for the hypergeometric series $_2\phi_1$. Our proof instead proceeds from the observation that the coefficients of  the ``homogeneous part'' of \eqref{E:Hqdiff} (which includes only those terms that have a factor $H(zq^k)$ for some $k\in\N$) have only even powers of $z$. This property allows us to solve the $q$-difference equation directly, even with the presence of the ``non-homogeneous'' final summation.

Writing the series expansion as $H(z) = \sum_{n \geq 0} \delta_n z^n$, \eqref{E:Hqdiff} implies that the coefficients $\delta_n$ satisfy the recurrence (for $n \geq 2$)
\begin{equation}
\label{E:deltarec}
\delta_n = \frac{\left(1 + t^{-1}q^{3n-4}\right)\left(1 + tq^{3n-2}\right)}{\left(1 - q^{3n-3}\right)\left(1 - q^{3n}\right)} \delta_{n-2}
+ \frac{(1 - \alpha - \beta)t (-1)^n q^{3n-2}}{\left(q^3; q^3\right)_n}.
\end{equation}
The initial values $\delta_0$ and $\delta_1$ are determined by using \eqref{E:Fdef} and \eqref{E:Hdef} to directly calculate
\begin{align*}
\sum_{n \geq 0} \delta_n z^n &= \gamma_0 + \left(\gamma_1 - \frac{\gamma_0}{1-q^3}\right) z + O\left(z^2\right)  = \beta + \frac{1+\alpha t q - \beta}{1-q^3} z + O\left(z^2\right).
\end{align*}
In particular, we read off
\begin{equation*}
\delta_0 = \beta \quad \text{and} \quad \delta_1 = \frac{1+\alpha t q - \beta}{1-q^3}.
\end{equation*}
Plugging in to \eqref{E:deltarec}, we find that for $n \geq 0$,
\begin{align}
\label{E:delta0}
\delta_{2n} & = \frac{\left(-t^{-1}q^2, -tq^4; q^6\right)_n}{\left(q^3; q^3\right)_{2n}}
\left(\beta + t(1-\alpha - \beta)\sum_{\ell = 1}^n \frac{q^{6\ell-2}}{\left(-t^{-1}q^2, -tq^4; q^6\right)_\ell}\right), \\
\label{E:delta1}
\delta_{2n+1} & = \frac{\left(-t^{-1}q^5, -tq^7; q^6\right)_n}{\left(q^3; q^3\right)_{2n+1}}
\left(1 -\beta + \alpha tq - t(1-\alpha - \beta)\sum_{\ell = 1}^n \frac{q^{6\ell+1}}{\left(-t^{-1}q^5, -tq^7; q^6\right)_\ell}\right).
\end{align}

For ease of calculation, we separate these two cases, setting
\begin{equation*}
H(z) = H_0(z) + H_1(z) := \sum_{n \geq 0} \delta_{2n} z^{2n} + \sum_{n \geq 0} \delta_{2n+1} z^{2n+1},
\end{equation*}
and similarly
\begin{equation*}
F_j(z) := \left(-z; q^3\right)_\infty H_j(z),
\end{equation*}
for $j = 0,1.$ We now plug \eqref{E:delta0} and \eqref{E:delta1} in to the definition of the $F_j$ and further expand the product $(-z; q^3)_\infty$ using \eqref{E:Euler2}. Collecting like powers of $z$, we obtain the series
\begin{align}
\label{E:F0}
F_0(z) & = \sum_{n \geq 0} \frac{z^n}{\left(q^3; q^3\right)_n}
\sum_{\substack{j, \, r \geq 0 \\ j + 2r = n}} q^{\frac{3j(j-1)}{2}}\qbinom{n}{j}_{q^3} \left(-t^{-1}q^2, -tq^4; q^6\right)_r \\
& \qquad \qquad \qquad \times \left(\beta + t(1 - \alpha - \beta)\sum_{\ell=1}^r \frac{q^{6\ell-2}}{\left(-t^{-1}q^2; -tq^4; q^6\right)_\ell}\right), \notag \\
\label{E:F1}
F_1(z) & = \sum_{n \geq 0} \frac{z^n}{\left(q^3; q^3\right)_n}
\sum_{\substack{j, \, r \geq 0 \\ j + 2r + 1 = n}} q^{\frac{3j(j-1)}{2}}\qbinom{n}{j}_{q^3} \left(-t^{-1}q^5, -tq^7; q^6\right)_r \\
& \qquad \times \left(1 - \beta + \alpha tq - t(1 - \alpha - \beta)\sum_{\ell=1}^r \frac{q^{6\ell+1}}{\left(-t^{-1}q^5; -tq^7; q^6\right)_\ell}\right). \notag
\end{align}
Recalling \eqref{E:gammadef} and \eqref{E:Fdef} and isolating the coefficients of $z^n$ in the above expressions completes the proof of Lemma \ref{L:Mainabfinite}.

\subsection{Proof of Theorem \ref{T:Mainab}}
\label{S:proof:infinite}

We now prove Theorem \ref{T:Mainab} by taking the limit as $n \to \infty$ of the expressions from Lemma \ref{L:Mainabfinite}.  We again separate the expressions arising from $H_0$ and $H_1$, letting $C_{0,3n-2}(t;q)$ and $C_{1,3n-2}(t;q)$ denote the inner sums in \eqref{E:F0} and \eqref{E:F1}, respectively. For example, we have
\begin{equation*}
F_0(z) = \sum_{n \geq 0} \frac{z^n}{\left(q^3; q^3\right)_n} \cdot C_{0, 3n-2}(t;q).
\end{equation*}

Using \eqref{E:qbinomLim}, the limit of these evaluates to
\begin{align}
C_0(t;q) & := \lim_{\substack{n \to \infty \\ n \text{ even}}} C_{0, 3n-2}(t;q)
\label{E:C0} \\
& = \left(-t^{-1}q^2, -tq^4; q^6\right)_\infty \sum_{\substack{j \geq 0 \\ j \text{ even}}} \frac{q^{\frac{3j(j-1)}{2}}}{\left(q^3; q^3\right)_j}
\left(\beta + t(1 - \alpha - \beta) \sum_{\ell \geq 1} \frac{q^{6\ell - 2}}{\left(-t^{-1}q^2, -tq^4; q^6\right)_\ell}\right). \notag
\end{align}

We are allowed to restrict to even $n$ because we know a priori (through combinatorial arguments) that the limit exists, and therefore has the same value as $n$ approaches $\infty$ along any subsequence. Using \eqref{E:CauchyEuler1}, the sum on $j$ evaluates to $(-q^3; q^3)_\infty$ (this evaluation was also noted in (4.8) of \cite{AAG95}).

For the final sum in \eqref{E:C0}, we first shift the summation index by $\ell \mapsto \ell + 1$ and then apply \eqref{E:Ramid} with $q \mapsto q^6, a = tq^4,$ and $b= t^{-1}q^2,$ obtaining
\begin{align}
\label{E:C0Ram}
& \sum_{\ell \geq 1} \frac{q^{6\ell - 2}}{\left(-t^{-1}q^2, -tq^4; q^6\right)_\ell}
= \frac{q^4}{\left(1 + t^{-1}q^2\right)\left(1 + tq^4\right)}
\sum_{\ell \geq 0} \frac{q^{6 \ell}}{\left(-t^{-1}q^8, -tq^{10}; q^6\right)_\ell} \\
& = \frac{q^4}{\left(1 + t^{-1}q^2\right)\left(1 + tq^4\right)}\left[
\left(1 + t^{-1}q^{-4}\right) \sum_{k \geq 0} \frac{(-1)^k q^{3k^2 + k} t^{-2k}}{\left(-t^{-1}q^8; q^6\right)_k} - \frac{t^{-1}q^{-4} \sum_{k \geq 0} (-1)^k q^{3k^2 +k} t^{-2k}}{\left(-t^{-1}q^8, -tq^{10}; q^6\right)_\infty} \right] \notag \\
& = t^{-1} \sum_{k \geq 0} \frac{(-1)^k q^{3k^2+k}t^{-2k}}{\left(-t^{-1}q^2; q^6\right)_{k+1}}
- \frac{t^{-1} \sum_{k \geq 0} (-1)^k q^{3k^2 + k}t^{-2k}}{\left(-t^{-1}q^2, -tq^4; q^6\right)_\infty}. \notag
\end{align}
To simplify further, consider the first sum on the third line of \eqref{E:C0Ram}. By shifting $k \mapsto k - 1$ and then applying \eqref{E:Rogid} with $q \mapsto q^6$ and $y = -t^{-1} q^{-4}$ we calculate
\begin{align}
\label{E:C0Rog}
\sum_{k \geq 0} \frac{(-1)^k q^{3k^2+k}t^{-2k}}{\left(-t^{-1}q^2; q^6\right)_{k+1}}
& = - \sum_{k \geq 1} \frac{(-1)^k q^{3k^2 - 5k + 2} t^{-2k + 2}}{\left(-t^{-1}q^2; q^6\right)_k} \\
& = -t^2 q^2 \left(-1 + \sum_{k \geq 0} t^{-3k} q^{9k^2 - 9k}\left(1 - t^{-2} q^{12k-2}\right)\right) \notag \\
& = \sum_{k\geq 0}\left(t^{-3k} q^{3k(3k+1)}-t^{-(3k+1)} q^{(3k+1)(3k+2)}\right). \notag
\end{align}
Combining \eqref{E:C0}, \eqref{E:C0Ram}, and \eqref{E:C0Rog}, we obtain the overall expression
\begin{align}
& \frac{C_0(t;q)}{\left(-q^3; q^3\right)_\infty} =
- (1 - \alpha - \beta) \sum_{k \geq 0} (-1)^k t^{-2k} q^{3k^2 + k} \label{E:C0final} \\
& + \left(-t^{-1}q^2, -tq^4; q^6\right)_\infty
\left(\beta + (1-\alpha-\beta)\sum_{k \geq 0} \left(t^{-3k}q^{3k(3k+1)} - t^{-(3k+1)}q^{(3k+1)(3k+2)}\right)\right). \notag
\end{align}

To complete the proof, we similarly evaluate the limit of the second summand from Lemma \ref{L:Mainabfinite} and simplify the resulting expressions. Proceeding as above, \eqref{E:qbinomLim} once more implies
\begin{align}
C_1(t;q) & := \lim_{\substack{n \to \infty \\ n \text{ odd}}} C_{1, 3n-2}(t;q)
= \left(-t^{-1}q^5, -tq^7; q^6\right)_\infty \sum_{\substack{j \geq 0 \\ j \text{ even}}} \frac{q^{\frac{3j(j-1)}{2}}}{\left(q^3; q^3\right)_j}
\label{E:C1} \\
& \qquad \qquad \times \left(1-\beta+\alpha t q - t(1 - \alpha - \beta) \sum_{\ell \geq 1} \frac{q^{6\ell + 1}}{\left(-t^{-1}q^5, -tq^7; q^6\right)_\ell}\right). \notag
\end{align}
The sum on $j$ is again evaluated by \eqref{E:CauchyEuler1}, and for the sum on $\ell$ we apply \eqref{E:Ramid} with $q \mapsto q^6, a = t^{-1} q^{-1},$ and $b = tq$, yielding
\begin{align}
\label{E:C1Ram}
\sum_{\ell \geq 1} \frac{q^{6\ell +1}}{\left(-t^{-1}q^5, -tq^7; q^6\right)_\ell}
& = q \left(-1 + \sum_{\ell \geq 0} \frac{q^{6 \ell}}{\left(-t^{-1}q^5, -tq^7; q^6\right)_\ell}\right) \\
& = -q + q(1 + tq) \sum_{k \geq 0} \frac{(-1)^k q^{3k^2+5k}t^{2k}}{\left(-tq^7; q^6\right)_k}
- \frac{tq^2 \sum_{k \geq 0} (-1)^k q^{3k^2 + 5k}t^{2k}}{\left(-t^{-1}q^5, -tq^7; q^6\right)_\infty}. \notag
\end{align}
Furthermore, applying \eqref{E:Rogid} with $q \mapsto q^6, y = -tq$ gives
\begin{equation}
\label{E:C1Rog}
\sum_{k \geq 0} \frac{(-1)^k q^{3k^2+5k}t^{2k}}{\left(-tq^7; q^6\right)_k}
= \sum_{k \geq 0} t^{3k} q^{9k^2 + 6k}\left(1 - t^{2} q^{12k+8}\right).
\end{equation}
Combining \eqref{E:C1}, \eqref{E:C1Ram}, and \eqref{E:C1Rog} results in the overall expression
\begin{align}
& \frac{C_1(t;q)}{\left(-q^3; q^3\right)_\infty} = (1 - \alpha - \beta) \sum_{k \geq 0} (-1)^k t^{2k+2} q^{3k^2 + 5k + 2} \label{E:C1final} \\
& + \left(-tq, -t^{-1}q^5;q^6\right)_\infty
\left(1-\beta - (1-\alpha-\beta)\sum_{k \geq 0} \left(t^{3k+1}q^{(3k+1)^2} - t^{3k+3}q^{(3k+3)^2}\right)\right). \notag
\end{align}
In order to write the inner sum in terms of $\Theta_2(t;q)$, note that $1 - \beta = \alpha + (1 - \alpha - \beta).$

The proof of Theorem \ref{T:Mainab} is complete once we add \eqref{E:C0final} and \eqref{E:C1final}. The final simplification comes from combining the first sum in both equations, using the fact that
\begin{equation*}
\sum_{k \geq 0} (-1)^k t^{2k+2} q^{3k^2 + 5k + 2} = - \sum_{k \leq -1} (-1)^k t^{-2k} q^{3k^2 + k}.
\end{equation*}
Recalling \eqref{E:theta}, this gives the theorem statement.

\end{document}